\documentclass[12pt,a4paper,reqno]{amsart}

\usepackage{anysize}
\usepackage[utf8]{inputenc}
\usepackage{amsmath}
\usepackage{amssymb}
\usepackage{array}

\newtheorem{theorem}{Theorem}[section]

\begin{document}

\title{WZ pairs and $q$-Analogues of Ramanujan series for $1/\pi$}

\author{Jesús Guillera}

\address{Department of Mathematics, University of Zaragoza, 50009 Zaragoza, SPAIN}

\email{jguillera@gmail.com}

\address{Department of Mathematics, IMAPP, Radboud University, PO Box 9010, 6500~GL Nijmegen, Netherlands}
\email{w.zudilin@math.ru.nl}

\date{}



\maketitle

\begin{abstract}
We prove q-analogues of two Ramanujan-type series for $1/\pi$ from $q$-analogues of ordinary WZ pairs.
\end{abstract}

\section{WZ-pairs}

A WZ pair (W=Wilf and Z=Zeilberger) is a pair formed by two hypergeometric (in their two variables) functions $F(n,k)$ and $G(n,k)$ such that
\[
F(n+1, k) -F(n, k) = G(n, k+1) - G(n, k).
\]
WZ pairs are gems because of their interesting properties and multiple applications. The property that we will use here is the following one: If $(F, G)$ is a WZ-pair and $F(0,k)=0$ and $F(\infty, k)=0$, for $k=0,1,2,\dots$, then
\[
\sum_{n=0}^{\infty} G(n,k) = C, \qquad k=0,1,2,\dots,
\]
where $C$ is a constant. Wilf and Zeilberger show that it was an immediate consequence of the definition if we sum for $n \geq 0$.
From a WZ-pair $(F_1, G_1)$ we can construct other $(F_2, G_2)$ pairs. We will use the following transformations:
\begin{equation}\label{pattern1}
F_2(n,k)=F_1(n,k+n), \quad G_2(n,k)=F_1(n+1,k+n)+G_1(n,k+n), 
\end{equation}
defined first in \cite{Guillera-generators},
\begin{equation}\label{pattern2}
F_2(n,k)=F_1(n,k-n), \quad G_2(n,k)=-F_1(n+1,k-n-1)+G_1(n,k-n),
\end{equation}
which is a backward version of (\ref{pattern1}), and 
\begin{equation}\label{pattern3}
F_2(n,k)=F_1(2n,k), \quad G_2(n,k)=G_1(2n,k)+G_1(2n+1,k),
\end{equation}
due to T. Amdeherhan and D. Zeilberger, which is the particular case $s=2$, $t=1$ of \cite[Formula 3]{Amde-Zeil}. Our strategy will consist in using the following result
\begin{equation}\label{our-strategy}
\sum_{n=0}^{\infty} G_1(n,k) = C, \quad \sum_{n=0}^{\infty} G_2(n,k) = C, \quad k=0,1,2, \dots,
\end{equation}
which holds when $F_1(0,k)=0$ if we assume the vanishing conditions at infinity (observe that then $F_2(0,k)=F_1(0,k)=0$ as well, and that in view of \cite[Theorem 7]{Zeil} the equality of the two constants arises). This strategy was first used in \cite{Guillera-generators} with ordinary WZ-pairs, and recently by Guo and Zudilin  in \cite{Guo-Zudilin} with $q$-WZ pairs using the iteration defined by (\ref{pattern1}). Of course we can combine two transformations and also iterate the process to get chains of WZ-pairs: $(F_3,G_3)$, $(F_4, G_4)$, etc. 

\section{$q$-analogues of Ramanujan series for $1/\pi$}

From a known $q$-analogue of a Ramanujan-type series for $1/\pi$, we obtain a $q$-analogue of another one.

\subsection{A known $q$-analogue of a Ramanujan series}

One of the WZ-pairs I found is \cite[Chapter 1, Pair 3.2]{Guillera-thesis}:
\begin{align*}
F(n,k) &= \frac{(-1)^n}{2^{3n+k}} \frac{(\frac12 -k)_n(\frac12+k)_n^2}{(1)_n^2(1+k)_n} \frac{(\frac12)_k}{(1)_k} \frac{16n^2}{2n-2k-1},
\\
G(n,k) &= \frac{(-1)^n}{2^{3n+k}} \frac{(\frac12 -k)_n(\frac12+k)_n^2}{(1)_n^2(1+k)_n} \frac{(\frac12)_k}{(1)_k} (6n+2k+1).
\end{align*}
Zudilin observed in \cite{Zu} that for proving supercongruences it is better to write it in the form
\[
F(n,k) = 8 \frac{(-1)^{n+k}}{2^{3n+k}} \frac{(\frac12)_{n-k-1}(\frac12)_{n+k}^2}{(1)_{n-1}^2(1)_{n+k}}, \quad 
G(n,k) = \frac{(-1)^{n+k}}{2^{3n+k}} \frac{(\frac12)_{n-k}(\frac12)_{n+k}^2}{(1)_n^2(1)_{n+k}} (6n+2k+1),
\]
and Guo and Liu proved in \cite{Guo-Liu} that
\begin{align*}
F(n,k) &= (-1)^{n+k} \frac{(q; q^2)_{n-k-1} (q; q^2)_{n+k}^2}{(q^4; q^4)_{n-1}^2 (q^4; q^4)_{n+k}} \frac{1}{1-q}, \\
G(n,k) &= (-1)^{n+k} \frac{(q; q^2)_{n-k} (q; q^2)_{n+k}^2}{(q^4; q^4)_{n}^2 (q^4; q^4)_{n+k}} [6n+2k+1],
\end{align*}
is a q-analogue of it. Replacing $q$ with $q^{-1}$ they wrote it in the more suitable form
\begin{align}
F(n,k) &= (-1)^{n+k} q^{(n+k)(3n-k)} \frac{(q; q^2)_{n-k-1} (q; q^2)_{n+k}^2}{(q^4; q^4)_{n-1}^2 (q^4; q^4)_{n+k}} \frac{1}{1-q}, \label{F-Guo}
\\ 
G(n,k) &= (-1)^{n+k} q^{(n+k)(3n-k)} \frac{(q; q^2)_{n-k} (q; q^2)_{n+k}^2}{(q^4; q^4)_{n}^2 (q^4; q^4)_{n+k}} [6n+2k+1], \label{G-Guo}
\end{align}
from which they derived the $q$-identity:
\begin{equation}\label{rama1-q-analog}
\sum_{n=0}^{\infty} (-1)^n q^{3n^2} [6n+1] \frac{(q;q^2)_n^3}{(q^4;q^4)_n^3} = \frac{(q^3;q^4)_{\infty}(q^5;q^4)_{\infty}}{(q^4;q^4)_{\infty}^2}, \qquad |q|<1.
\end{equation}
Writing it in the following way:
\[
\sum_{n=0}^{\infty} (-1)^n q^{3n^2} [6n+1] \frac{(q;q^2)_n^3}{(q^4;q^4)_n^3} = \lim_{n \to \infty} \frac{(q^3;q^4)_{n}(q^5;q^4)_{n}}{(q^4;q^4)_{n}^2}, \qquad |q|<1,
\]
and using the property
\[
\lim_{q \to 1^{-}} \frac{(q^a;q^b)_n}{(1-q)^n} = b^n \left( \frac{a}{b} \right)_n,
\]
one can easily check that (\ref{rama1-q-analog}) is a q-analogue of the Ramanujan-type series of level $4$:
\begin{equation}
\sum_{n=0}^{\infty} \frac{(\frac12)_n^3}{(1)_n^3}(6n+1) \frac{(-1)^n}{2^{3n}} = \frac{2\sqrt 2}{\pi}.
\end{equation}
We have used standard notation \cite[Introduction]{Guo-Zudilin}.  Guo and Zudilin \cite{Guo-Zudilin} used the iteration pattern (\ref{pattern1}) to derive $q$-analogues of some other Ramanujan and Ramanujan-type series. 

\subsection{A $q$-analogue of a Ramanujan series of level $1$}

We add a new $q$-identity to the list:
\begin{multline}\label{new-q-analogue}
\sum_{n=0}^{\infty} (-1)^n q^{7n^2} \frac{(q; q^2)_{3n}(q; q^2)^2_n}{(q^4; q^4)_{2n}^2 (q^4; q^4)_n} \left( \frac{[6n+1]}{[4n+4]} \frac{[6n+3]q^{14n+7}-[10n+7]q^{10n+3}}{(1+q^{2n+1})^2 (1+q^{4n+2})^2} + [10n+1] \right) \\ = \frac{(q^3;q^4)_{\infty}(q^5;q^4)_{\infty}}{(q^4;q^4)_{\infty}^2},
\end{multline}
which is a $q$-analogue of the Ramanujan-type series for $1/\pi$ of level one:
\begin{equation}\label{rama-level-1}
\sum_{n=0}^{\infty} \frac{(\frac12)_n(\frac16)_n(\frac56)_n}{(1)_n^3} (154n+15) \left( - \frac38 \right)^{3n} = \frac{32 \sqrt 2}{\pi}.
\end{equation}
\begin{proof}
Combining the transformations (\ref{pattern2}) and (\ref{pattern3}), we get the following WZ pair:
\begin{align}
F_2(n,k) &= F_1(2n, k-n), \label{F2-combine} \\
G_2(n,k) &= -F_1(2n+2, k-n-1) + G_1(2n, k-n) + G_1(2n+1, k-n). \label{G2-combine}
\end{align}
If $(F_1, G_1)$ is the WZ pair in (\ref{F-Guo}, \ref{G-Guo}), we get
\begin{equation}\label{G1-eq-G2}
\sum_{n=0}^{\infty} G_2(n,k) = \sum_{n=0}^{\infty} G_1(n,k) =  \frac{(q^3;q^4)_{\infty}(q^5;q^4)_{\infty}}{(q^4;q^4)_{\infty}^2}, \quad k=0,1,2\dots.
\end{equation}
From (\ref{G1-eq-G2}) and
\[
\sum_{n=0}^{\infty} G_2(n,0) = \sum_{n=0}^{\infty} \left[ -F_1(2n+2,-n-1) +G_1(2n, -n) + G_1(2n+1, -n) \right],
\]
we obtain (\ref{new-q-analogue}). Finally, taking the limit as $q \to 1^{-}$ we see that (\ref{new-q-analogue}) is a $q$-analogue of the identity (\ref{rama-level-1}).
\end{proof}

\subsection{Generalization with an extra parameter $a$}
We can generalize Guo's $q$-WZ-pair (\ref{F-Guo}, \ref{G-Guo}), with an extra parameter $a$, in the following way:
\begin{align*}
F_1(n,k) &= (-1)^{n+k} q^{(n+k)(3n-k)} \frac{(q; q^2)_{n-k-1} (aq; q^2)_{n+k} (q/a; q^2)_{n+k}}{(q^4; q^4)_{n+k}(aq^4; q^4)_{n-1} (q^4/a;q^4)_{n-1}} \frac{1}{1-q}, 
\\ 
G_1(n,k) &= (-1)^{n+k} q^{(n+k)(3n-k)} \frac{(q; q^2)_{n-k} (aq; q^2)_{n+k}(q/a; q^2)_{n+k}}{(q^4; q^4)_{n+k} (aq^4; q^4)_n (q^4/a; q^4)_{n}} [6n+2k+1].
\end{align*}
Then, from the last formula of the proof of \cite[Theorem 4.4]{Guo-Zudilin-2}:
\[
\sum_{n=0}^{\infty} (-1)^{n} q^{3n^2} \frac{(q; q^2)_{n} (aq; q^2)_{n}(q/a; q^2)_{n}}{(q^4; q^4)_{n} (aq^4; q^4)_n (q^4/a; q^4)_{n}} [6n+1]=\frac{(q^3;q^4)_{\infty}(q^5;q^4)_{\infty}}{(aq^4;q^4)_{\infty} (q^4/a; q^4)_{\infty}},
\]
and using the WZ-pair in (\ref{F2-combine}, \ref{G2-combine}), we obtain the following generalization of (\ref{new-q-analogue}):
\begin{multline}\label{q-analogue-1-a}
\frac{(q^3;q^4)_{\infty}(q^5;q^4)_{\infty}}{(aq^4;q^4)_{\infty}(q^4/a;q^4)_{\infty}}=\sum_{n=0}^{\infty} (-1)^n q^{7n^2} \frac{(q; q^2)_{3n}(aq; q^2)_n(q/a; q^2)_n}{(q^4; q^4)_{n} (aq^4; q^4)_{2n} (q^4/a; q^4)_{2n}} \times \\ 
\left\{ \frac{(1-aq^{2n+1})(1-q^{2n+1}/a)}{(1-aq^{8n+4})(1-q^{8n+4}/a)} \frac{[6n+1]}{[4n+4]} \left( [6n+3]q^{14n+7}-[10n+7]q^{10n+3}\right) + [10n+1] \right\},
\end{multline}
which has an extra parameter $a$ as those used in \cite{Guo-Zudilin-2} to prove supercongruences.

\section{$q$-analogues of other Ramanujan series for $1/\pi$}

From a known $q$-analogue of a Ramanujan series for $1/\pi$, we obtain a $q$-analogue of another one.

\subsection{A known $q$-analogue of a Ramanujan series}
From a $q$-hypergeometric identity, Guo and Zudilin \cite{Guo-Zudilin} found the following $q$-analogue
\begin{equation}\label{q-guo-zud}
\sum_{n=0}^{\infty} q^{2n^2}  \frac{(q;q^2)_n^2(q;q^2)_{2n}}{(q^6;q^6)_n^2(q^2;q^2)_{2n}} [8n+1] = \frac{(q^3;q^2)_{\infty}(q^3;q^6)_{\infty}}{(q^2;q^2)_{\infty}(q^6;q^6)_{\infty}},
\end{equation}
of the Ramanujan series of level $2$:
\[
\sum_{n=0}^{\infty}\frac{\left(\frac12\right)_n\left(\frac14\right)_n\left(\frac34\right)_n}{(1)_n^3}(8n+1) \frac{1}{9^n} = \frac{2 \sqrt 3}{\pi}.
\]

\subsection{A $q$-analogue of another Ramanujan series}
We can rewrite the WZ pair $7$ of \cite[p. 28]{Guillera-thesis} in the following way:
\begin{align}
F(n,k) &= 18 (-1)^k \frac{\left(\frac12\right)_{n+k}^2\left(\frac12\right)_{2n-k-1} \left(\frac12\right)_k}{(1)_{n-1}^2(1)_{2n+2k}} \frac{3^k}{9^{n}}, \label{F-8n1} \\
G(n,k) &= (-1)^k \frac{\left(\frac12\right)_{n+k}^2\left(\frac12\right)_{2n-k}\left(\frac12\right)_k}{(1)_{n}^2(1)_{2n+2k}} \frac{3^k}{9^{n}} (8n+2k+1) \label{G-8n1},
\end{align}
Inspired by (\ref{q-guo-zud}), and after some trials for the exponent of $q^{2n^2+?}$, we get the following $q$-analogue of the WZ-pair (\ref{F-8n1},\ref{G-8n1}):
\begin{align}
F_1(n,k) &= \frac{(-1)^k}{1-q} q^{2n^2+4nk-k^2} \frac{(q;q^2)_{n+k}^2(q;q^2)_{2n-k-1}}{(q^6;q^6)_{n-1}^2(q^2;q^2)_{2n+2k}} 
(q^3;q^6)_k , \label{F1-q-8n1} \\
G_1(n,k) &= (-1)^k q^{2n^2+4nk-k^2} \frac{(q;q^2)_{n+k}^2(q;q^2)_{2n-k}}{(q^6;q^6)_{n}^2(q^2;q^2)_{2n+2k}} (q^3;q^6)_k [8n+2k+1], \label{G1-q-8n1}
\end{align}
In addition, we see that $F_1(0,k)=0$. Then, for the pair ($F_2(n,k), G_2(n,k)$) defined in (\ref{pattern1}), we have
\[
\sum_{n=0}^{\infty} G_2(n,0) = \sum_{n=0}^{\infty} G_2(n,k) = \sum_{n=0}^{\infty} G_1(n,k) = \sum_{n=0}^{\infty} G_1(n,0),
\]
and it leads to
\begin{multline}
\sum_{n=0}^{\infty} (-1)^n q^{5n^2} \frac{(q;q^2)_{2n}^2 (q;q^2)_n(q^3;q^6)_n}{(q^6;q^6)_n^2(q^2;q^2)_{4n}} 
\left( \frac{q^{8n+2}[4n+1]}{(1+q^{2n+1})(1+q^{4n+1})(1+q^{4n+2})} + [10n+1] \right) \\
=\frac{(q^3;q^2)_{\infty}(q^3;q^6)_{\infty}}{(q^2;q^2)_{\infty}(q^6;q^6)_{\infty}},
\end{multline}
which is a q-analogue of the Ramanujan series of level $2$:
\[
\sum_{n=0}^{\infty}\frac{\left(\frac12\right)_n\left(\frac14\right)_n\left(\frac34\right)_n}{(1)_n^3}(28n+3) \left( \frac{-1}{48} \right)^n = \frac{16 \sqrt 3}{3\pi}.
\]

\subsection{Generalization with an extra parameter $a$}

We can generalize the WZ pair (\ref{F1-q-8n1},\ref{G1-q-8n1}), in the following way
\begin{align}
F_1(n,k) &= \frac{(-1)^k}{1-q} q^{2n^2+4nk-k^2} \frac{(aq;q^2)_{n+k} (q/a;q^2)_{n+k}(q;q^2)_{2n-k-1}}{(aq^6;q^6)_{n-1}(q^6/a;q^6)_{n-1}
(q^2;q^2)_{2n+2k}} (q^3;q^6)_k , \label{F1-q-8n1-a} \\
G_1(n,k) &= (-1)^k q^{2n^2+4nk-k^2} \frac{(aq;q^2)_{n+k} (q/a;q^2)_{n+k} (q;q^2)_{2n-k}}{(aq^6;q^6)_{n} (q^6/a;q^6)_n (q^2;q^2)_{2n+2k}} (q^3;q^6)_k [8n+2k+1]. \label{G1-q-8n1-a}
\end{align}
Then, if ($F_2(n,k), G_2(n,k)$) is the WZ-pair defined in (\ref{pattern1}), we have 
\[
\sum_{n=0}^{\infty} G_2(n,0)=\sum_{n=0}^{\infty} G_2(n,k) = \sum_{n=0}^{\infty} G_1(n,k) = \sum_{n=0}^{\infty} G_1(n,0).
\]
But Guo and Zudilin proved in \cite{Guo-Zudilin-2} that
\[
\sum_{n=0}^{\infty} G_1(n,0) = \frac{(q^5;q^6)_{\infty}(q^7;q^6)_{\infty}(aq^3;q^6)_{\infty}(q^3/a;q^6)_{\infty}}{(q^2;q^6)_{\infty}(q^4;q^6)_{\infty}(aq^6;q^6)_{\infty}(q^6/a;q^6)_{\infty}},
\]
and we obtain
\begin{multline}\label{q-analogue-28n3-a}
 \sum_{n=0}^{\infty} (-1)^n q^{5n^2} \frac{(aq;q^2)_{2n} (q/a;q^2)_{2n} (q;q^2)_n(q^3;q^6)_n}{(aq^6;q^6)_n (q^6/a;q^6)_n (q^2;q^2)_{4n}}  \times \\ \left( \frac{(1-aq^{4n+1})(1-q^{4n+1}/a)}{(1-q)^2} \frac{q^{8n+2}}{[8n+2](1+q^{2n+1})(1+q^{4n+2})} + [10n+1] \right) \\ =
 \frac{(q^5;q^6)_{\infty}(q^7;q^6)_{\infty}(aq^3;q^6)_{\infty}(q^3/a;q^6)_{\infty}}{(q^2;q^6)_{\infty}(q^4;q^6)_{\infty}(aq^6;q^6)_{\infty}(q^6/a;q^6)_{\infty}},
\end{multline}
which has an extra parameter $a$ as those used in \cite{Guo-Zudilin-2} to prove supercongruences.

\section*{Acknowledgement}
I thank Wadim Zudilin for encourage me to introduce the extra parameter $a$ that one needs to prove the kind of supercongruences in \cite{Guo-Zudilin-2}.

\newpage

\setcounter{section}{0}
\def\thesection{\Alph{section}}

\section{{\large Appendix. ($q$)-Congruences}}

\begin{center}
\bf Wadim Zudilin
\end{center}

\vskip 1cm

\subsection{Congruences from (\ref{q-analogue-1-a})}
Let $c_q(a;n)$ be the term inside the summation of (\ref{q-analogue-1-a}).
Take any odd $m>0$. Let $d>1$ be a divisor of $m$, and let $\zeta=\zeta_d$ be a primitive $d$-th root of unity.
Then the right-hand side in \eqref{q-analogue-1-a} tends to 0 as $q$ tends to $\zeta$ radially, because of the numerator
$(q^3;q^2)_\infty=\prod_{j=1}^\infty(1-q^{2j+1})$. At the same time
$$
\lim_{q\to\zeta}c_q(a;n)=c_\zeta(a;n),
$$
and the latter vanishes if $n>\lfloor(d-1)/6\rfloor$. Therefore, the limit as $q\to\zeta$ in \eqref{q-analogue-1-a} results in
$$
\sum_{n=0}^{\lfloor(d-1)/6\rfloor}c_\zeta(a;n)=0
$$
implying, in particular, that
$$
\sum_{n=0}^{(m-1)/2}c_\zeta(a;n)=0
\quad\text{and}\quad
\sum_{n=0}^{m-1}c_\zeta(a;n)=0.
$$
The equalities mean that
$$
\sum_{n=0}^{(m-1)/2}c_q(a;n)\equiv0\pmod{\Phi_d(q)}
\quad\text{and}\quad
\sum_{n=0}^{m-1}c_q(a;n)\equiv0\pmod{\Phi_d(q)},
$$
where $\Phi_d(q)$ is the $d$-th cyclotomic polynomial.
Because the congruences are true for any divisor $d$ of $m$, we conclude that
\begin{equation}\label{cyclo-cong1}
\sum_{n=0}^{(m-1)/2}c_q(a;n)\equiv0\pmod{[m]}
\quad\text{and}\quad
\sum_{n=0}^{m-1}c_q(a;n)\equiv0\pmod{[m]}.
\end{equation}

Now, if we substitute $a=q^m$ or $a=q^{-m}$ in \eqref{q-analogue-1-a}, then the hypergeometric sum terminates at $n=(m-1)/2$
and we obtain the equalities
\begin{align*}
&
\sum_{n=0}^{(m-1)/2}c_q(q^m;n)
=\sum_{n=0}^{(m-1)/2}c_q(q^{-m};n)
=\sum_{n=0}^{m-1}c_q(q^m;n)
=\sum_{n=0}^{m-1}c_q(q^{-m};n)
\\ &\quad
=\frac{(q^3;q^2)_{\infty}}{(q^{4+m};q^4)_{\infty}(q^{4-m};q^4)_{\infty}}
=(-q)^{(m-1)(m-3)/8}[m].
\end{align*}
This evaluation means that
\begin{equation}\label{cong1-zud}
\sum_{n=0}^{(m-1)/2}c_q(a;n)\equiv(-q)^{(m-1)(m-3)/8}[m]
\quad\text{and}\quad
\sum_{n=0}^{m-1}c_q(a;n)\equiv(-q)^{(m-1)(m-3)/8}[m]
\end{equation}
modulo $(a-q^m)(1-aq^m)$. Combining the result with the congruences \eqref{cyclo-cong1} we deduce that the congruences
\eqref{cong1-zud} are valid modulo $[m](a-q^m)(1-aq^m)$. Finally, considering the limiting case as $a\to1$ and arguing as in \cite{Guo-Zudilin-2} we arrive at the following result for the partial sums of
\begin{align*}
&
c_q(1;n)=(-1)^nq^{7n^2}
\frac{(q;q^2)_{3n}(q;q^2)_n^2}{(q^4;q^4)_n(q^4;q^4)_{2n}^2}
\\ &\;
\times
\bigg(\frac{[2n+1]^2[6n+1]}{[4n+4]\,[8n+4]^2}\,
\big([6n+3]q^{14n+7}-[10n+7]q^{10n+3}\big)+[10n+1]\bigg).
\end{align*}

\begin{theorem} \label{th:cong1}
For any $m>1$ odd, we have
\begin{align*}
\sum_{n=0}^{(m-1)/2}c_q(1;n)&\equiv(-q)^{(m-1)(m-3)/8}[m]\pmod{[m]\Phi_m(q)^2}
\\ \intertext{and}
\sum_{n=0}^{m-1}c_q(1;n)&\equiv(-q)^{(m-1)(m-3)/8}[m]\pmod{[m]\Phi_m(q)^2}.
\end{align*}
In particular, as $q\to1^-$ we obtain the Ramanujan-type supercongruences
\begin{align*}
\sum_{n=0}^{(p-1)/2}\frac{(\frac12)_n(\frac16)_n(\frac56)_n}{n!^3}(154n+15)\biggl(-\frac38\biggr)^{3n}
&\equiv15p\,\biggl(\frac{-2}p\biggr)\pmod{p^3}
\\ \intertext{and}
\sum_{n=0}^{p-1}\frac{(\frac12)_n(\frac16)_n(\frac56)_n}{n!^3}(154n+15)\biggl(-\frac38\biggr)^{3n}
&\equiv15p\,\biggl(\frac{-2}p\biggr)\pmod{p^3}
\end{align*}
valid for any prime $p>2$.
\end{theorem}

In the last passage we have also used $(-1)^{(m-1)(m-3)/8}=\bigl(\frac{-2}m\bigr)$ for odd $m$.

\subsection{Congruences from (\ref{q-analogue-28n3-a})}
The situation for the identity (\ref{q-analogue-28n3-a}) is somewhat different from the one considered above and in \cite{Guo-Zudilin-2}. Let $c_q(a;n)$ be the term inside the summation. The presence of $(q;q^2)_n(q^3;q^6)_n/(q^2;q^2)_{4n}$ makes the radial limits of individual terms $c_q(a;n)$ as $q\to\zeta$ in the sum infinite, starting from a certain index, for any root of unity $\zeta$.
However, we still have terminating sums at $a=q^m$ and $a=q^{-m}$ for any $m>0$ coprime with 6,
and this allows us to conclude that
\begin{equation}\label{cong1}
\sum_{n=0}^{(m-1)/2}c_q(a;n)\equiv q^{-(m-1)/2}[m]\biggl(\frac{-3}m\biggr)
\quad\text{and}\quad
\sum_{n=0}^{m-1}c_q(a;n)\equiv q^{-(m-1)/2}[m]\biggl(\frac{-3}m\biggr)
\end{equation}
modulo $(a-q^m)(1-aq^m)$ for all such $m$ (see \cite[Lemma 3.1]{Guo-Zudilin-2} for evaluation of the right-hand sides).
By performing the limit as $a\to1$ we therefore arrive at a weaker version of the expected congruences for the partial sums of
\begin{align*}
c_q(1;n)
&=(-1)^nq^{5n^2}
\frac{(q;q^2)_n(q^3;q^6)_n(q;q^2)_{2n}^2}{(q^2;q^2)_{4n}(q^6;q^6)_n^2}
\\ &\quad\times
\bigg(\frac{[4n+1]^2q^{8n+2}}{(1+q^{2n+1})(1+q^{4n+2})[8n+2]}+[10n+1]\bigg).
\end{align*}

\begin{theorem} \label{th:cong2}
For any $m>1$ relatively prime with $6$, we have
\begin{align*}
\sum_{n=0}^{(m-1)/2}c_q(1;n)&\equiv q^{-(m-1)/2}[m]\biggl(\frac{-3}m\biggr)\pmod{\Phi_m(q)^2}
\\ \intertext{and}
\sum_{n=0}^{m-1}c_q(1;n)&\equiv q^{-(m-1)/2}[m]\biggl(\frac{-3}m\biggr)\pmod{\Phi_m(q)^2}.
\end{align*}
In particular, as $q\to1^-$ we obtain the Ramanujan-type supercongruences
\begin{align*}
\sum_{n=0}^{(p-1)/2}\frac{(\frac12)_n(\frac14)_n(\frac34)_n}{n!^3}(28n+3)\biggl(-\frac1{48}\biggr)^n
&\equiv3p\,\biggl(\frac{-3}p\biggr)\pmod{p^2}
\\ \intertext{and}
\sum_{n=0}^{p-1}\frac{(\frac12)_n(\frac14)_n(\frac34)_n}{n!^3}(28n+3)\biggl(-\frac1{48}\biggr)^n
&\equiv3p\,\biggl(\frac{-3}p\biggr)\pmod{p^2}
\end{align*}
valid for any prime $p>3$.
\end{theorem}


\begin{thebibliography}{99}

\bibitem{Amde-Zeil}
T. Amdeberhan and D. Zeilberger,
Hypergeometric acceleration via the WZ-method,
The Electronic Journal of Combinatorics {\bf 4} (1997).

\bibitem{Guillera-generators}
J. Guillera, 
Generators of some Ramanujan formulas,
Ramanujan J. {\bf 11}, 2006, 41–48.

\bibitem{Guillera-thesis}
J. Guillera,
Series de Ramanujan: Generalizaciones y conjeturas [PhD thesis].
Zaragoza (Spain): Universidad de Zaragoza; 2007.

\bibitem{Guo-Liu}
V. Guo and J-C Liu,
$q$-Analogues of two Ramanujan-type formulas for $1/\pi$, 
J. Diff. Equ. Appl. (to appear), 5 pages,
Preprint: https://arxiv.org/abs/1802.01944.

\bibitem{Guo-Zudilin}
V. Guo and W. Zudilin,
Ramanujan-type formulas for $1/\pi$: $q$-analogues,
Integral Transforms Spec. Functions (to appear), 9 pages,  
Preprint: https://arxiv.org/abs/1802.04616.

\bibitem{Guo-Zudilin-2}
V. Guo and W. Zudilin,
A $q$-microscope for supercongruences,
Preprint [March 2018, 24 pages] at https://arxiv.org/abs/1803.01830.

\bibitem{Zeil} 
D. Zeilberger,
Closed form (pun intended!),
Contemp. Math. {\bf 143}, (1993), 579--608.

\bibitem{Zu}
W. Zudilin,
Ramanujan-type supercongruences,
J. Number Theory 129:8 (2009), 1848–1857.

\end{thebibliography}
\end{document}